\documentclass[12pt]{article}

\usepackage{graphicx}  
\usepackage{amsmath}  
\usepackage{amssymb,amsthm}
\usepackage{color}
\usepackage{tikz}
\usepackage{pgfplots}

\newcommand{\R}{{\mathbb R}}\newcommand{\N}{{\mathbb N}}

\let\epsilon\varepsilon
\let\theta\vartheta
\let\hat\widehat

\newtheorem{theorem}{Theorem}[section]

\newtheorem{remark}[theorem]{Remark}

\title{A note on the validity of the
Schr\"odinger approximation for the Helmholtz equation}
\author{Maximilian Klumpp, Guido Schneider \\
{\small
Institut f\"ur Analysis, Dynamik und Modellierung,} \\ {\small Universit\"at Stuttgart, Pfaffenwaldring 57, } \\ {\small 70569 Stuttgart, Germany}}

\begin{document}

\maketitle

\begin{abstract}
Time-harmonic electromagnetic waves in vacuum
are described by the Helmholtz equation
$
\Delta u+\omega ^{2}u=0 
$ for $ (x,y,z) \in \R^3 $.
For the evolution of such waves along
the $z$-axis a Schr\"odinger equation can be
derived through a multiple scaling
ansatz. It is the purpose of this paper
to justify this formal approximation by
proving bounds between this formal approximation
and true solutions of the original system.
The challenge of the presented validity
analysis is the fact that the Helmholtz
equation is ill-posed as an evolutionary
system along the $z$-axis.
\end{abstract}


\section{Introduction}

Electromagnetic waves are described by the system
of Maxwell's equations \cite{Jackson}. In vacuum 
this system of equations reduce
to the scalar linear wave equation
\[
\partial ^{2}_{t}u=\Delta u=\partial ^{2}_{x}u+\partial ^{2}_{y}u+\partial ^{2}_{z}u ,
\]
with
$u\left(x,y,z,t\right) ,
x,y,z,t \in \mathbb{R}  $.
For time harmonic waves
\[
u\left( x,y,z,t\right) =v\left( x,y,z\right) e^{i\omega t} 
\]
we find
\begin{equation}
	\label{eq2}
-\omega ^{2}v=\partial ^{2}_{x}v+\partial ^{2}_{y}v+\partial ^{2}_{z}v .
\end{equation}
We are interested in the evolution of waves along
the z-axis and so we write  \eqref{eq2}
as evolutionary
system w.r.t. the $z$-variable.
The resulting system
\begin{equation}
	\label{eq3}
	\partial ^{2}_{z}v=-\omega ^{2}v-\partial ^{2}_{x}v-\partial ^{2}_{y}v 
\end{equation}
is solved by harmonic waves
\[
v\left( x,y,z\right) =e^{\lambda z+ik_{x}x+ik_{y}y},
\]
with $
\lambda(k_x,k_y), k_x,k_y\in\R$ and 
\begin{equation} \label{lambda}
\lambda^{2}=k^{2}_{x}+k^{2}_{y}-\omega ^{2} .
\end{equation}
 Since \eqref{eq2}
is an
elliptic system the $z$-evolutionary system \eqref{eq3}
in general is an ill-posed initial value problem.
Out of this ill-posed system a well-posed initial
value problem can be derived through some
multiple scaling analysis. In detail, inserting the
ansatz
\[
v\left( x,y,z\right) \approx \psi _{app}\left( x,y,z\right)  =e^{ik_{z}z}w\left( X,Y,Z\right)   ,
\]
with 
\begin{equation}\label{scaling}
X=\varepsilon x,\quad Y=\varepsilon y, \quad Z=\varepsilon ^{2}z,
\end{equation} 
and $ 0 < \varepsilon \ll 1 $ a small perturbation parameter,
yields
\[
-k^{2}_{z}w+2i\varepsilon ^{2}k_{z}\partial _{Z}w+\varepsilon ^{4}\partial ^{2}_{Z}w =
-\omega ^{2}w-\varepsilon ^{2} \partial ^{2}_{X}w-\varepsilon ^{2} \partial ^{2}_{Y}w ,
\]
see the subsequent Remark \ref{rem3} for an additional discussion.
Choosing $\omega ^{2}=k^{2}_{z}$
and ignoring the terms of
order $\mathcal{O}\left( \varepsilon ^{4}\right)  $ gives that $w $ has to satisfy in
lowest order the Schr\"odinger equation
\begin{equation}\label{eq4}
	2ik_{z}\partial _{Z}w=-\partial ^{2}_{X}w-\partial ^{2}_{Y}w ,
\end{equation}
which, in contrast to the Helmholtz equation \eqref{eq3}, is a well-posed initial value problem w.r.t. the evolutionary variable $Z $.
This formal approach is widely used in the physics and engineering literature for instance to describe the evolution of light beams in vacuum and media \cite{schmidt,rafferty},  to compute beam quality factors 
in step-index fibers \cite{yoda},   or in the modeling of physical optics phenomena by complex
ray tracing \cite{harvey}. It is known as the "paraxial approximation".  A motivation of this approximation is given in  \cite{Marte}, but, to our knowledge, no rigorous analysis has been made so far  to justify this formal approximation.
Therefore, it is the purpose of this note to prove rigorous bounds between this formal approximation and true solutions of the original system for small $ \varepsilon > 0 $. 
Such estimates are known for dispersive/hyperbolic systems already for many years, see Remark 
\ref{rem23}.

\section{The approximation result}

{\bf Notation.} Many possibly different constants are denoted
with the symbol $C $ if they can be chosen
independently of the small perturbation parameter
$0 < \varepsilon \ll 1 $.
The Fourier transform of a function
$u $ w.r.t. the variables $x $ and $y$ is denoted with $\widehat {u} $
and is given by
\[
\widehat {u}\left( k_{x},k_{y}\right) =\dfrac {1}{\left( 2\pi \right) ^{2}}\int_{\mathbb{R}^2} u\left( x,y\right)  
e^{-ik_{x}x-ik_{y}y}d\left( x,y\right)  .
\]
The Sobolev space $H^s$ is equipped with the
norm
\[
\left\| u\right\| _{H^{s}}=\sum_{m_{1}+m_{2}\leq s,m_{1},m_{2}\in \N_{0} }\left( \int_{\mathbb{R}^2} \left| \partial ^{m_{1}}_{x}\partial ^{m_{2}}_{y}u\right| ^{2}d\left( x,y\right) \right) ^{1/2}.
\]
We introduce the weighted Lebesgue space $L^2_s$
equipped with the norm
\[
\left\| \widehat {u}\right\| _{L^{2}_s}=\left( \int_{\mathbb{R}^2} \left| \widehat {u}\left( k_{x},k_{y}\right) \right| ^{2}\rho\left( \left| k\right| \right) ^{2s}d\left( k_{x},k_{y}\right) \right) ^{1/2} ,
\]
where $\rho\left( \left| k\right| \right) =\left( 1+\left| k\right| ^{2}\right) ^{1/2} 
$ 
and $\left| k\right| ^{2}=k^{2}_{x}+k^{2}_{y} $.
We use that Fourier transform is an isomorphism
between $H^s$ and $L^2_s$, i.e., for $ s \in \mathbb{N}_0 $ there exist $C_1,C_2>0$
such that for all $u \in H^s$
\begin{align}
\left\| \widehat {u}\right\| _{L^{2}_{s}}\leq C_{1}\left\| u\right\| _{H^{s}}\leq C_{2}\left\| \widehat {u}\right\| _{L^{2}_{s}} .\label{equiv}
\end{align}
See \cite[Lemma 7:3.31]{SU17}. Hence, for all $ s \geq 0 $ we redefine $ \| u \|_{H^s}  = \| \widehat{u} \|_{L^2_s}$. 

Then we prove the following approximation result.
\begin{theorem}\label{thhaupt}
	Fix $s_{A}\geq \max \left( 4,s\right)$ and $Z_0>0$. Let
	${w\in C\left( \left[ 0,Z_{0}\right] ,H^{s_A}\right)}$
be a solution of the Schr\"odinger equation \eqref{eq4}.
Then there exist a $C> 0 $ and an
$\varepsilon _{0} > 0 $
such that for all 
$
\varepsilon \in \left( 0,\varepsilon _{0}\right)  $
there are solutions
$v $ of  \eqref{eq3}
with
\[
\sup_{z\in \left[ 0,Z_{0}/\varepsilon ^{2}\right] } \left\| v\left( x,y,z\right) -e^{ik_{z}z}w\left( \varepsilon x,\varepsilon y,\varepsilon ^{2}z\right) \right\| _{H^{s}\left( dx,dy\right) }\leq C\varepsilon  ,
\]
In particular, by Sobolev's embedding theorem we have
\[
\sup_{z\in \left[ 0,Z_{0}/\varepsilon ^{2}\right] } 
\sup_{x,y \in \R^2}
\left| v\left( x,y,z\right) -e^{ik_{z}z}w\left( \varepsilon x,\varepsilon y,\varepsilon ^{2}z\right) \right|\leq C\varepsilon  .
\]
\end{theorem}
\begin{remark}{\rm 
	The challenge of the presented validity
analysis is the fact that in every $ H^s $ the Helmholtz
equation is ill-posed as an evolutionary
system along the $z$-axis.}
\end{remark}
\begin{remark} \label{rem23}{\rm 
	The method presented in this note does not apply to nonlinear
problems since it uses a cut-off in Fourier space
which is not respected by nonlinear terms.
Hence, it cannot be used to justify the NLS
equation for weakly nonlinear elliptic Maxwell
models. See for instance \cite{Ka88,KSM92,Schn05}
for validity results for the NLS approximation
of nonlinear dispersive/hyperbolic systems and the textbooks \cite{Rauch,SU17} for an introduction and a recent overview.}
\end{remark}
\begin{remark} \label{rem3}{\rm 
Alternatively to the scaling (\ref{scaling}) we can introduce the small  parameter 
$ 0 < \varepsilon \ll 1 $
already at the beginning by considering
\[
u\left( x,y,z,t\right) =v\left( x,y,z\right) e^{i\omega t/\varepsilon^2} 
\]
which can be physically motivated by the highly oscillatory character of light.
The
ansatz is then given by 
\[
 \psi _{app}\left( x,y,z\right)  =e^{ik_{z}z/\varepsilon^2}w\left( x,y,z\right)   ,
\]
which
yields
as above the Schr\"odinger equation \eqref{eq4}. It is an easy exercise to reformulate
Theorem \ref{thhaupt} and the subsequent proof w.r.t. this scaling.
}
\end{remark}

\section{The proof}

We have that for fixed $ \omega $ the eigenvalues $ \lambda $ defined by \eqref{lambda}
are purely imaginary for $ k^{2}_{x}+k^{2}_{y}\leq \omega ^{2} $ and real-valued 
for $ k^{2}_{x}+k^{2}_{y}\geq \omega ^{2} $
with $ |\lambda(k_x,k_y) | \to \infty $ for $ k^{2}_{x}+k^{2}_{y}  \to \infty $.
Since \eqref{eq3} is then an ill-posed initial value problem
we use a cut-off function in Fourier space
to remove the ill-posed part of \eqref{eq3}. In detail,
we define a projection $P_{hyp} $ on the hyperbolic part
of \eqref{eq3}  by
\begin{align}
P_{hyp}u=\mathcal{F}^{-1}\chi
\mathcal{F}u \label{P}
\end{align}
where 
\[
\chi(k)=\left\{\begin{array}{cl}
1 ,& \text{for}	\quad k^{2}_{x}+k^{2}_{y}\leq \omega ^{2}/2 ,\\
0 , & \text{else}.
\end{array}\right.
\]
Moreover, let $P_{ell}=I-P_{hyp} $. 
For all $ s \geq 0 $  we have 
$$ 
\| P_{hyp} u \|_{H^s} \leq \| \chi \widehat{u} \|_{L^2_s} \leq \| \widehat{u} \|_{L^2_s} \leq \|  u \|_{H^s}
$$
and similarly $ \| P_{ell} u \|_{H^s} \leq \|  u \|_{H^s}$.
For the subsequent
estimates it turns out to be advantageous to take $\omega^2/2$
instead of $\omega^2$ in the definition of $P_{hyp}$. In order to estimate
the difference between solutions $v $ of \eqref{eq3}
and the Schr\"odinger approximation 
$\psi_{app}$ we
work with  the hyperbolic part $P_{hyp} \psi_{app}$ of the
Schr\"odinger approximation. Then we use the triangle
inequality to estimate
\[
\left\| v-\psi_{app}\right\| _{H^{s}}\leq \left\| v-P_{hyp}\psi _{app}\right\| _{H^{s}}+\left\| P_{hyp}\psi _{app}-\psi _{app}\right\| _{H^{s}} .
\]
The two terms on the right hand side are estimated
in the following two subsections.

\subsection{Estimating the well-posed part}

\label{sec4}

We consider  \eqref{eq3} with initial conditions $P_{hyp}\psi _{app}|_{z=0}$ and $\partial_zP_{hyp}\psi _{app}|_{z=0}$.
Since \eqref{eq3} is a linear system, the Fourier
support of the solution $v $ of \eqref{eq3}
is preserved.
The difference $R = v -P_{hyp} \psi_{app}$ then satisfies
\begin{equation} \label{eqR}
\partial_z^2 R = - \partial_x^2 R-\partial_y^2 R- \omega^2 R - \varepsilon^4 e^{ik_zz}P_{hyp}\partial_Z^2  w,
\end{equation}
with initial data $R|_{z=0}=0$ and $\partial_zR|_{z=0}=0$.
 The inhomogeneity
can be expressed via the right hand side
of the Schr\"odinger equation, namely
\[
\partial_Z^2  w= -\partial_Z  \left(\frac{1}{2ik_z}(\partial_X^2 w+\partial_Y^2 w)\right)=-
\frac{1}{4k_z^2}(\partial_X^2 +\partial_Y^2 )^2 w.
\]
For 
$ w \in C([0,Z_0],H^{s_A}) $  
we thus have
$ \partial_Z^2  w\in C([0,Z_0],H^{s_A-4}) $,  
i.e., 
 there exist $C_1, C_{res}> 0 $
such that
\[
\| \varepsilon^4 P_{hyp}\partial_Z^2  w \|_{L^2} \leq C_1\varepsilon^4\|(\partial_X^2+\partial_Y^2)^2w\|_{L^2}\le C_{res}\varepsilon^4 \| w \|_{H^{4}}\le C_{res}\varepsilon^4 \| w \|_{H^{s_A}}
\]
for all $ z \geq 0 $.
Since for
functions with compact Fourier support each
$ H^s $-norm 
 can be estimated by the
$ L^2 $-norm, in particular for $ u $ with $ u = P_{hyp} u  $ we have 
\begin{equation} \label{newl}
\|  u \|_{H^s} = \| \chi \widehat{u} \|_{L^2_s} \leq \| \chi \widehat{u} \rho^s  \|_{L^2_0}
\leq  \| \chi  \rho^s  \|_{L^{\infty}} \| \widehat{u}   \|_{L^2_0} \leq C \| u \|_{L^2},
\end{equation}
with a constant $ C_s $
for each fixed
$s $.
Therefore, we finally
have
$$ 
\| \varepsilon^4 P_{hyp}\partial_Z^2  w \|_{H^s}  \leq 
C_s \| \varepsilon^4 P_{hyp}\partial_Z^2  w \|_{L^2} 
\leq C_s C_{res}\varepsilon^4 \| w \|_{H^{s_A}}.
$$
In Fourier space, w.r.t. $ x $ and $y $,
\eqref{eqR}
 is given by
\begin{equation} \label{eqB}
\partial_z^2 \widehat{R} = -  \widehat{\omega}^2 \widehat{R} - \varepsilon^4e^{ik_zz} \chi \partial_Z^2 \widehat{w} , 
\end{equation}
with $\widehat{\omega}^2(k_x,k_y) =  \omega ^{2} -k^{2}_{x}-k^{2}_{y}$ and $\hat{R}, \hat{w}$ is the Fourier transform w.r.t. the variables $X, Y$.
Multiplying
\eqref{eqB}
 with 
$ \overline{\partial_z \widehat{R}} $
and integrating w.r.t.
$ k_z $ and $ k_y $ yields
\begin{align*}
\frac{d}{dz} \int_{\mathbb{R}^2} | \partial_z \widehat{R} |^2 + | \widehat{\omega} \widehat{R} |^2 
d(k_x,k_y) & \leq   2|\textrm{Re} \int_{\mathbb{R}^2} \overline{\partial_z \widehat{R} }(\varepsilon^4 \chi \partial_Z^2 \widehat{w} )d(k_x,k_y) |\\
& \leq  \varepsilon^2 \int_{\mathbb{R}^2} | \partial_z \widehat{R} |^2d(k_x,k_y)
+  \varepsilon^6 \int_{\mathbb{R}^2} | \chi \partial_Z^2 \widehat{w} |^2 d(k_x,k_y) .
\end{align*}
The second summand can be estimated by
\begin{align*}
\varepsilon^6 \int_{\mathbb{R}^2} | \chi \partial_Z^2 \widehat{w} |^2 d(k_x,k_y)&=\varepsilon^6\|\chi\partial_Z^2 \widehat{w}\|_{L^2(dk_x,dk_y)}^2=\varepsilon^6\|P_{hyp}\partial_z^2w\|_{L^2(dx,dy)}^2\\
&\le\varepsilon^4\|P_{hyp}\partial_z^2w\|_{L^2(dX,dY)}^2\le C_{res}\varepsilon^4\|w\|_{H^{s_A}(dX,dY)}.
\end{align*}

Thus, for
$ E = \int | \partial_z \widehat{R} |^2 + | \widehat{\omega} \widehat{R} |^2 
d(k_x,k_y) $ 
 we find
\[
\frac{d}{dz} E \leq \varepsilon^2 E  +  \varepsilon^4 C_{Res}^2 \| w \|^2_{ C([0,Z_0],H^{4})}.
\]
Gronwall's inequality then yields
\[
E(z)  \leq \int_0^z e^{\varepsilon^2 (z-s)} \varepsilon^4 C_{Res}^2 \| w \|^2_{ C([0,Z_0],H^{4})} ds
\le C \varepsilon^2 \| w \|^2_{ C([0,Z_0],H^{4})}
\]
for $ z \in [0,Z_0/\varepsilon^2] $.
Since by construction $ \widehat{\omega} \chi$ is bounded away from zero,
 $E(z)^{1/2}$  is an upper bound for
the $L^2$-norm of $ \widehat{R} $ and thus by Parseval's
inequality  an upper bound for the $L^2$-norm of $ R $, in detail
\begin{align*}
	\|R\|_{L^2}=\|\hat{R}\|_{L^2}\le C_3\|\hat{\omega}\hat{R}\|_{L^2}\le C_3E(z)^{1/2}
\end{align*}
for some $C_3>0$. Since $ R $ has
a compact support in Fourier space, as in \eqref{newl}  the
$L^2$-estimate implies an $ H^s $-estimate for $ R $ for every $ s \geq 0 $. 

\subsection{Estimating the ill-posed part}
\label{sec5}

For estimating the term $ \left\| P_{hyp}\psi _{app}-\psi _{app}\right\| _{H^{s}}  $ we use
that 
$$
\widehat{\psi} _{app}(k_x,k_y,z)= \varepsilon^{-2} \widehat{w}\left(\frac{k_x}{\varepsilon},\frac{k_y}{\varepsilon},
  \varepsilon^2 z\right) e^{ik_{z}z}
$$
is strongly concentrated at the wave vector $ (k_x,k_y)=(0,0) $. 
For $ s_A \geq s $ 
we estimate
\begin{align*}
\lefteqn{\left\| P_{hyp}\psi _{app}-\psi _{app}\right\| _{H^{s}} }\\  &\overset{(\ref{equiv})}{\le} C 
 \left\| \widehat{P}_{hyp}\widehat{\psi} _{app}- \widehat{\psi}_{app}\right\| _{L^2_s}
 \overset{(\ref{P})}{=} C \left\| (1-\chi)  \widehat{\psi}_{app}\right\| _{L^2_s}  \\
 & =  C \left\|(1-\chi(k)) \varepsilon^{-2} \widehat{w}\left(\frac{k_x}{\varepsilon},\frac{k_y}{\varepsilon},
 \varepsilon^2 z\right){\rho(k)^s} \right\| _{L^2(dk_x,dk_y)} \\
 & \le  C \sup_{k \in \R^2} \left| (1-\chi(k)) \frac{\rho(k)^s}{\rho\left(\frac{k}{\varepsilon}\right)^{s_A}} \right|
  \left\| \varepsilon^{-2} \widehat{w}\left(\frac{k_x}{\varepsilon},\frac{k_y}{\varepsilon},
  \varepsilon^2 z\right) \rho\left(\frac{k}{\varepsilon}\right)^{s_A}
 \right\| _{L^2(dk_x,dk_y)} \\
 & \le  C \varepsilon^{s_A-1}   \| \widehat{w}(K_x,K_y) \|_{L^2_{s_A}(dK_x,dK_y)} 
  \overset{(\ref{equiv})}{\leq}  C \varepsilon^{s_A-1} \| w \|_{H^{s_A}},
\end{align*}
with $K_x=\frac{k_x}{\varepsilon}, K_y=\frac{k_y}{\varepsilon}$, where the loss of $  \varepsilon^{-1}$ comes from the scaling properties
of the $L^2$-norm and of Fourier transform.

\subsection{Summary}

Combining the estimates of Section \ref{sec4} and Section \ref{sec5}
gives via the triangle inequality
\begin{align*}
	\left\| v-\psi_{app}\right\| _{H^{s}}& \le  \left\| v-P_{hyp}\psi _{app}\right\| _{H^{s}}+\left\| P_{hyp}\psi _{app}-\psi _{app}\right\| _{H^{s}}  \\
	& \le C_3\sqrt{C \varepsilon^2 } \| w \|_{ C([0,Z_0],H^{4})}
	+  C \varepsilon^{s_A-1} \| w \|_{C([0,Z_0],H^{s_A})} \\
	& \le  \tilde{C}\varepsilon
\end{align*}
for $  z \in [0,Z_0/\varepsilon^2] $ and $s_{A}\geq \max \left( 4,s \right) $.
Therefore, we are done. \qed

%
%
%
%

\bibliographystyle{plain}
\bibliography{nlsellipticrev}

\end{document}